# Viscosity solutions for systems of parabolic variational inequalities

LUCIAN MATICIUC[1], ETIENNE PARDOUX[2], AUREL RĂŞCANU[3,4] and
ADRIAN ZĂLINESCU[4]

[1]*Department of Mathematics, Gheorge Asachi Technical University of Iaşi, Bd. Carol I, no. 11, Iaşi 700506, Romania. E-mail: lucianmaticiuc@yahoo.com*

[2]*LATP, Université de Provence, CMI, 39, rue F. Joliot-Curie, 13453, Marseille Cedex 13, France. E-mail: pardoux@cmi.univ-mrs.fr*

[3]*Faculty of Mathematics, Alexandru Ioan Cuza University of Iaşi, Bd. Carol I, no. 9, Iaşi 700506, Romania. E-mail: aurel.rascanu@uaic.ro*

[4]*Octav Mayer Institute of Mathematics, Iaşi, Bd. Carol I, no. 8, Iaşi 700506, Romania. E-mail: adrian.zalinescu@gmail.com*

In this paper, we first define the notion of viscosity solution for the following system of partial differential equations involving a subdifferential operator:

$$\begin{cases} \dfrac{\partial u}{\partial t}(t,x) + \mathcal{L}_t u(t,x) + f(t,x,u(t,x)) \in \partial\varphi(u(t,x)), & t\in[0,T), x\in\mathbb{R}^d, \\ u(T,x) = h(x), & x\in\mathbb{R}^d, \end{cases}$$

where $\partial\varphi$ is the subdifferential operator of the proper convex lower semicontinuous function $\varphi:\mathbb{R}^k \to (-\infty,+\infty]$ and $\mathcal{L}_t$ is a second differential operator given by $\mathcal{L}_t v_i(x) = \frac{1}{2}\operatorname{Tr}[\sigma(t,x)\sigma^*(t,x)\mathrm{D}^2 v_i(x)] + \langle b(t,x), \nabla v_i(x)\rangle$, $i\in\overline{1,k}$.

We prove the uniqueness of the viscosity solution and then, via a stochastic approach, prove the existence of a viscosity solution $u:[0,T]\times\mathbb{R}^d \to \mathbb{R}^k$ of the above parabolic variational inequality.

*Keywords:* Feynman–Kac formula; systems of variational inequalities; viscosity solutions

## 1. Introduction

Viscosity solutions theory was introduced by M.G. Crandall and P.L. Lions [3]. This theory allows one to study nonlinear equations which admit as solutions continuous functions, without any further smoothness constraints. The classical work in the field of viscosity solutions for second order partial differential equations is the survey paper of M.G. Crandall, H. Ishii and P.L. Lions [2], where the authors provide several equivalent ways to define the notion of viscosity solution, as well as some very general existence and uniqueness theorems. Starting with the study of backward stochastic differential equations (introduced by E. Pardoux and S. Peng [10]), generalized Feynman–Kac representation formulas have been obtained for the viscosity solutions of semilinear partial differential equations. In [4], R.W.R. Darling and E. Pardoux studied elliptic equations with Dirichlet boundary conditions; furthermore, in [7], Y. Hu treated elliptic equations with homogeneous Neumann boundary conditions. E. Pardoux and S. Zhang extended these results in [13] for the case of parabolic systems with nonlinear Neumann boundary conditions.

N. El Karoui *et al.* [5] considered the case of reflected solutions of one-dimensional backward stochastic differential equations, related to an obstacle problem for a parabolic partial differential equation. The more general case of backward stochastic differential equations involving subdifferential operators and the connection with parabolic variational inequalities has been studied by E. Pardoux and A. Răşcanu in [11] and [12].

The aim of this paper is to consider the more general case of variational inequalities for systems of partial differential equations. The paper is organized as follows. In Section 2, we define the extended notion of viscosity solution for a system of parabolic variational inequalities. We then formulate the existence and uniqueness result and prove uniqueness. In Section 4, we use a stochastic approach in order to prove the existence result; by using a backward stochastic variational inequality, we obtain a probabilistic formula for the viscosity solution of our system.

## 2. Main results

The goal of this paper is to study the existence and uniqueness of the viscosity solution of the following system of parabolic variational inequalities:

$$\begin{cases} \dfrac{\partial u}{\partial t}(t,x) + \mathcal{L}_t u(t,x) + f(t,x,u(t,x)) \in \partial\varphi(u(t,x)), & t \in [0,T), x \in \mathbb{R}^d, \\ u(T,x) = h(x), & x \in \mathbb{R}^d, \end{cases} \quad (1)$$

where $\mathcal{L}_t v$, with $v \in C^2(\mathbb{R}^d; \mathbb{R}^k)$, is given by

$$(\mathcal{L}_t v)_i(x) = \frac{1}{2}\operatorname{Tr}[\sigma(t,x)\sigma^*(t,x)\mathrm{D}^2 v_i(x)] + \langle b(t,x), \mathrm{D}v_i(x)\rangle$$

$$= \frac{1}{2}\sum_{j,l=1}^d (\sigma\sigma^*)_{jl}(t,x)\frac{\partial^2 v_i(x)}{\partial x_j \partial x_l} + \sum_{j=1}^d b_j(t,x)\frac{\partial v_i(x)}{\partial x_j}, \qquad i \in \overline{1,k},$$

and $T > 0$ is the fixed finite horizon.

We make the following standard assumptions.

(A.1) The functions $b\colon [0,T]\times\mathbb{R}^d \to \mathbb{R}^d$ and $\sigma\colon [0,T]\times\mathbb{R}^d \to \mathbb{R}^{d\times d}$ are Lipschitz with constant $L$.

(A.2) The function $f:[0,T]\times\mathbb{R}^d\times\mathbb{R}^k\to\mathbb{R}^k$ is continuous and there exists $\gamma\in\mathbb{R}$ such that

$$\langle y-\tilde{y}, f(t,x,y)-f(t,x,\tilde{y})\rangle \leq \gamma|y-\tilde{y}|^2$$

for all $t\in[0,T]$, $x\in\mathbb{R}^d$, $y,\tilde{y}\in\mathbb{R}^k$.

(A.3) The function $h:\mathbb{R}^d\to\mathbb{R}^k$ is continuous and there exist some $M_1>0$, $p\geq 0$ such that, for all $t\in[0,T]$, $x\in\mathbb{R}^d$,

$$|h(x)|\leq M_1(1+|x|^p).$$

(A.4) The function $\varphi:\mathbb{R}^k\to(-\infty,+\infty]$ is proper (that is, $\varphi\not\equiv +\infty$), convex, lower semicontinuous (l.s.c.) and there exist $M_2>0$ and $r\geq 0$ such that

$$|\varphi(h(x))|\leq M_2(1+|x|^r)\qquad \forall x\in\mathbb{R}^d.$$

We recall that the subdifferential $\partial\varphi$ is defined by

$$\partial\varphi(u)=\{u^*\in\mathbb{R}^k:\langle u^*,v-u\rangle+\varphi(u)\leq\varphi(v),\forall v\in\mathbb{R}^k\}.$$

It is a common practice to regard $\partial\varphi$ as a subset of $\mathbb{R}^k\times\mathbb{R}^k$, by writing $(u,u^*)\in\partial\varphi$ instead of $u^*\in\partial\varphi(u)$. We define

$$\mathrm{Dom}(\varphi)=\{u\in\mathbb{R}^k:\varphi(u)<+\infty\},$$
$$\mathrm{Dom}(\partial\varphi)=\{u\in\mathbb{R}^k:\partial\varphi(u)\neq\varnothing\}.$$

We now give a result which allows us to define the notion of a directional derivative of a convex function (see [1] for more details).

**Theorem 1.** *Let $\varphi:\mathbb{R}^k\to(-\infty,+\infty]$ be a convex function. Then, for all $u\in\mathrm{Dom}(\varphi)$ and $z\in\mathbb{R}^k$, there exist*

$$\varphi'_-(u;z):=\lim_{t\nearrow 0}\frac{\varphi(u+tz)-\varphi(u)}{t}=\sup_{t<0}\frac{\varphi(u+tz)-\varphi(u)}{t},$$
$$\varphi'_+(u;z):=\lim_{t\searrow 0}\frac{\varphi(u+tz)-\varphi(u)}{t}=\inf_{t>0}\frac{\varphi(u+tz)-\varphi(u)}{t}.\qquad(2)$$

*Moreover, the following hold:*

(a) $\varphi'_-(u;z)\leq\varphi'_+(u;z),\forall u\in\mathrm{Dom}(\varphi),\forall z\in\mathbb{R}^k$;
(b) $\varphi'_-(u;-z)=-\varphi'_+(u;z),\forall u\in\mathrm{Dom}(\varphi),\forall z\in\mathbb{R}^k$;
(c) $\varphi'_-(u;\cdot)$ is superlinear and $\varphi'_+(u;\cdot)$ is sublinear;
(d) if $u$ and $z$ are such that there exists $\delta>0$ such that $u+tz\in\mathrm{Dom}(\varphi),\forall t\in(-\delta,\delta)$, then $\varphi'_-(u;z),\varphi'_+(u;z)\in\mathbb{R}$.

If we take $k = 1$, then we know that, in every point $u \in \text{Dom}(\varphi)$,

$$\partial \varphi(u) = \mathbb{R} \cap [\varphi'_-(u), \varphi'_+(u)], \tag{3}$$

where $\varphi'_-(u)$ and $\varphi'_+(u)$ are, respectively, the left derivative and right derivative of $\varphi$ at the point $u$.

The following proposition generalizes the above characterization to the case of $k \geq 1$.

**Proposition 2.** *Let $\varphi : \mathbb{R}^k \to (-\infty, +\infty]$ be a proper convex function and $u \in \text{Dom}(\varphi)$. The following statements are equivalent:*

(i) $u^* \in \partial \varphi(u)$;
(ii) $\langle u^*, z \rangle \geq \varphi'_-(u; z), \forall z \in \mathbb{R}^k$;
(iii) $\langle u^*, z \rangle \leq \varphi'_+(u; z), \forall z \in \mathbb{R}^k$.

**Proof.** We first prove that (i) implies assertion (ii). Let $u^* \in \partial \varphi(u)$. From the definition, we have that

$$\langle u^*, y - u \rangle + \varphi(u) \leq \varphi(y) \qquad \forall y \in \mathbb{R}^k. \tag{4}$$

For all $t < 0$ and any direction $z \in \mathbb{R}^k$, we obtain, using (4) with $y = u + tz$, that

$$\langle u^*, tz \rangle + \varphi(u) \leq \varphi(u + tz) \qquad \forall z \in \mathbb{R}^k, \forall t < 0$$

$$\Leftrightarrow \quad \langle u^*, z \rangle \geq \frac{\varphi(u + tz) - \varphi(u)}{t} \qquad \forall z \in \mathbb{R}^k, \forall t < 0.$$

Passing to the limit as $t \to 0$, $t < 0$, we obtain that

$$\langle u^*, z \rangle \geq \varphi'_-(u; z) \qquad \forall z \in \mathbb{R}^k.$$

Suppose now that (ii) holds, that is,

$$\langle u^*, z \rangle \geq \varphi'_-(u; z) \qquad \forall z \in \mathbb{R}^k.$$

From (2), we deduce that

$$\langle u^*, z \rangle \geq \frac{\varphi(u + tz) - \varphi(u)}{t} \qquad \forall z \in \mathbb{R}^k, \forall t < 0.$$

Consequently

$$\langle u^*, tz \rangle \leq \varphi(u + tz) - \varphi(u) \qquad \forall z \in \mathbb{R}^k, \forall t < 0.$$

Hence,

$$\langle u^*, y - u \rangle \leq \varphi(y) - \varphi(u) \qquad \forall y \in \mathbb{R}^k,$$

that is, $u^* \in \partial \varphi(u)$.

The equivalence between (i) and (iii) follows in the same manner. □

Let us define, for $u \in \overline{\mathrm{Dom}(\varphi)}$ and $z \in \mathbb{R}^k$,

$$\varphi'_*(u;z) = \liminf_{\substack{v \to u \\ v \in \mathrm{Dom}(\partial\varphi)}} \varphi'_-(v;z), \qquad \varphi'^{,*}(u;z) = \limsup_{\substack{v \to u \\ v \in \mathrm{Dom}(\partial\varphi)}} \varphi'_+(v;z).$$

For $u \in \mathbb{R}^k$, let (with the usual convention that $\inf \varnothing = +\infty$)

$$|\partial\varphi|_0(u) = \inf |\partial\varphi(u)|.$$

If $u \in \mathrm{Dom}(\partial\varphi)$, then there exists a unique $u^* \in \mathbb{R}^k$, denoted $(\partial\varphi)_0(u)$, such that $|\partial\varphi|_0(u) = |(\partial\varphi)_0(u)|$.

We can now give, using Proposition 2, the definition of a viscosity solution for the multidimensional parabolic variational inequality (1). First, we denote the set of $d \times d$ symmetric real matrices by $S^d$ and give the definitions of the superjet and the subjet of a function.

**Definition 3.** *Let $\bar{u}:[0,T] \times \mathbb{R}^d \to \mathbb{R}$ be a continuous function and $(t,x) \in [0,T] \times \mathbb{R}^d$. We denote by $\mathcal{P}^{2,+}\bar{u}(t,x)$ (the* parabolic superjet *of $\bar{u}$ at $(t,x)$) the set of triples $(p,q,X) \in \mathbb{R} \times \mathbb{R}^d \times S^d$ which are such that for all $(s,y) \in [0,T] \times \mathbb{R}^d$ in a neighborhood of $(t,x)$,*

$$\bar{u}(s,y) \leq \bar{u}(t,x) + p(s-t) + \langle q, y-x \rangle$$
$$+ \tfrac{1}{2}\langle X(y-x), y-x\rangle + \mathrm{o}(|s-t| + |y-x|^2).$$

*We similarly define $\mathcal{P}^{2,-}\bar{u}(t,x)$ (the* parabolic subjet *of $\bar{u}$ at $(t,x)$) as the set of triples $(p,q,X) \in \mathbb{R} \times \mathbb{R}^d \times S^d$ which are such that for all $(s,y) \in [0,T] \times \mathbb{R}^d$ in a neighborhood of $(t,x)$,*

$$\bar{u}(s,y) \geq \bar{u}(t,x) + p(s-t) + \langle q, y-x \rangle$$
$$+ \tfrac{1}{2}\langle X(y-x), y-x\rangle + \mathrm{o}(|s-t| + |y-x|^2).$$

*Here, $r \mapsto \mathrm{o}(r)$ denotes any function such that $\lim_{r \to 0} \frac{\mathrm{o}(r)}{r} = 0$.*

**Definition 4.** *Let $u:[0,T] \times \mathbb{R}^d \to \mathbb{R}^k$ be a continuous function satisfying $u(T,\cdot) = h(\cdot)$. The function $u$ is a viscosity solution of (1) if*

$$u(t,x) \in \overline{\mathrm{Dom}(\varphi)} \qquad \forall (t,x) \in [0,T) \times \mathbb{R}^d$$

*and if, for all $z \in \mathbb{R}^k$, at any point $(t,x) \in [0,T) \times \mathbb{R}^d$, for any $(p,q,X) \in \mathcal{P}^{2,+}\langle u,z \rangle(t,x)$, we have that*

$$p + \tfrac{1}{2}\mathrm{Tr}((\sigma\sigma^*)(t,x)X) + \langle b(t,x), q \rangle + \langle f(t,x,u(t,x)), z \rangle \geq \varphi'_*(u(t,x);z).$$

**Remark 5.** *If $u:[0,T] \times \mathbb{R}^d \to \mathbb{R}^k$ is a continuous function which satisfies $u(T,\cdot) = h(\cdot)$, then $u$ is a viscosity solution of (1) if and only if*

$$u(t,x) \in \overline{\mathrm{Dom}(\varphi)} \qquad \forall (t,x) \in [0,T) \times \mathbb{R}^d$$

and if, for all $z \in \mathbb{R}^k$, at any point $(t,x) \in [0,T) \times \mathbb{R}^d$, for any $(p,q,X) \in \mathcal{P}^{2,-}\langle u,z\rangle(t,x)$, we have that

$$p + \tfrac{1}{2}\operatorname{Tr}((\sigma\sigma^*)(t,x)X) + \langle b(t,x), q\rangle + \langle f(t,x,u(t,x)), z\rangle \leq \varphi'^{,*}(u(t,x); z).$$

We now present the main results.

**Theorem 6 (Existence).** *Let assumptions* (A.1)–(A.4) *be satisfied. The multidimensional parabolic variational inequality (1) then has at least a viscosity solution.*

For the proof of the existence theorem, we use a stochastic approach; we study a backward stochastic variational inequality which allows us to give a probabilistic formula for the solutions of (1). We present this approach in the last section.

Concerning the uniqueness for equation (1), we need to impose more restrictive assumptions.

(A.5) For all $u \in \operatorname{Dom}(\varphi)$, there exists a neighborhood $V$ of $u$ such that $(\partial\varphi)_0$ is bounded on $\operatorname{Dom}(\partial\varphi) \cap V$.

(A.6) If $u \in \operatorname{Dom}(\varphi)$ and $z \in \mathbb{R}^k$ such that $u + z \in \operatorname{Dom}(\varphi)$, then there exists a neighbourhood $V$ of $u$ such that

$$\forall v \in V \cap \operatorname{Dom}(\partial\varphi),\ \exists t > 0: \qquad v + tz \in \operatorname{Dom}(\partial\varphi).$$

**Theorem 7 (Uniqueness).** *Let the assumptions* (A.1), (A.2), (A.5) *and* (A.6) *be satisfied. The multidimensional parabolic variational inequality (1) then has at most one viscosity solution in the class of continuous functions with polynomial growth.*

## 3. Proof of the uniqueness theorem

We first prove the following lemma.

**Lemma 8.** *If $\varphi \colon \mathbb{R}^k \to (-\infty, +\infty]$ is a convex function, proper and l.s.c., then we have:*

(a) *for every $u, v \in \operatorname{Dom}(\varphi)$,*

$$\varphi'_-(u, u-v) \geq \varphi'_+(v, u-v); \tag{5}$$

(b) *for every $u, v \in \operatorname{Dom}(\partial\varphi)$ and $z \in \mathbb{R}^k$, $\varphi'_*(u;z) \leq \varphi'_-(u;z)$ and $\varphi'^{,*}(u;z) \geq \varphi'_+(u;z)$;*

(c) *for every $z \in \mathbb{R}^k$, $\varphi'_*(\cdot;z)$ is l.s.c. on $\operatorname{Dom}(\varphi)$ and $\varphi'^{,*}(\cdot;z)$ is upper semicontinous (u.s.c.) on $\operatorname{Dom}(\varphi)$;*

(d) *for every $u, v \in \overline{\operatorname{Dom}(\varphi)}$, if*

$$\limsup_{\substack{(u',v') \to (u,v) \\ u',v' \in \operatorname{Dom}(\partial\varphi)}} \left(\inf_{t>0}[|\partial\varphi|_0(u' - t(u-v)) + |\partial\varphi|_0(v' + t(u-v))]\right) < +\infty, \tag{6}$$

*then*

$$\varphi'_*(u, u - v) \geq \varphi'^{,*}(v, u - v).$$

**Proof.** (a) By the convexity of $\varphi$, for $t \in (0,1)$, we have

$$\varphi((1-t)u + tv) + \varphi(tu + (1-t)v) \leq \varphi(u) + \varphi(v).$$

Therefore,

$$-\frac{\varphi(u - t(u-v)) - \varphi(u)}{t} \geq \frac{\varphi(v + t(u-v)) - \varphi(v)}{t}$$

and passing to the limit as $t \searrow 0$, we obtain $\varphi'_-(u, u-v) \geq \varphi'_+(v, u-v)$.

Results (b) and (c) follow in a standard way from the definitions of $\varphi'_*$ and $\varphi'^{,*}$, so we skip their proofs.

(d) For simplicity, let us define $z = u - v$. From (6), we can find $M > 0$ and $r > 0$ such that, for $u' \in B(u, r) \cap \mathrm{Dom}(\partial\varphi)$, $v' \in B(v, r) \cap \mathrm{Dom}(\partial\varphi)$, there is a $t_{u',v'} < 0$ such that, for $t \in (t_{u',v'}, 0)$, we have

$$|\partial\varphi|_0(u' + tz) + |\partial\varphi|_0(v' - tz) \leq M.$$

Let $u'$, $v'$ and $t$ be as above. Then,

$$\frac{\varphi(u' + tz) - \varphi(u')}{t} = \frac{\varphi(u' + t(u' - v')) - \varphi(u')}{t} + \frac{\varphi(u' + tz) - \varphi(u' + t(u' - v'))}{t}$$

$$\geq \frac{\varphi(u' + t(u' - v')) - \varphi(u')}{t} + \langle (\partial\varphi)_0(u' + tz), z - (u' - v') \rangle.$$

It follows that

$$\frac{\varphi(u' + tz) - \varphi(u')}{t} \geq \frac{\varphi(u' + t(u' - v')) - \varphi(u')}{t} - M|z - (u' - v')|.$$

Passing to the limit as $t \nearrow 0$, we have

$$\varphi'_-(u'; z) \geq \varphi'_-(u'; u' - v') - M|z - (u' - v')|.$$

In a similar manner, we obtain that

$$\varphi'_+(v'; u' - v') \geq \varphi'_+(v'; z) - M|z - (u' - v')|.$$

Combining the last two inequalities with (5), we deduce that

$$\varphi'_-(u'; z) \geq \varphi'_+(v'; z) - 2M|z - (u' - v')|.$$

Passing to the limit as $u' \to u$ and $v' \to v$, we obtain that $\varphi'_*(u; z) \geq \varphi'^{,*}(v; z)$. $\square$

***Remark 9.*** From point (c) of the above lemma, we can consider in Definition 4, $\widetilde{\mathcal{P}}^{2,+}\langle u, z\rangle$ (resp., $\widetilde{\mathcal{P}}^{2,-}\langle v, z\rangle$) instead of $\mathcal{P}^{2,+}\langle u, z\rangle$ (resp., $\mathcal{P}^{2,-}\langle v, z\rangle$), where $\widetilde{\mathcal{P}}^{2,+}$ and $\widetilde{\mathcal{P}}^{2,-}$ denote the respective closures of operators $\mathcal{P}^{2,+}$ and $\mathcal{P}^{2,-}$ (see the definitions in [2]).

**Proof.** Let $u, v \in C([0,T] \times \mathbb{R}^d; \mathbb{R}^k)$ be viscosity solutions for (1). We must prove that $u = v$ on $[0,T] \times \mathbb{R}^d$.

For some $\lambda > 0$ to be made precise below, we carry out the following transformations:
$$\bar{u}(t,x) = e^{\lambda t}\xi^{-1}(x)u(t,x), \qquad \bar{v}(t,x) = e^{\lambda t}\xi^{-1}(x)v(t,x),$$

where $\xi(x) = (1+|x|^2)^{\mu/2}$, with $\mu$ larger than the order of growth of $u$ and $v$.

Then, $\bar{u}$ and $\bar{v}$ are bounded and they are solutions, in the viscosity sense, of

$$\begin{cases} \dfrac{\partial \bar{u}(t,x)}{\partial t} + \breve{\mathcal{E}}_t \bar{u}(t,x) + \bar{f}(t,x,\bar{u}(t,x)) \in e^{\lambda t}\xi^{-1}(x)\,\partial\varphi(e^{-\lambda t}\xi(x)\bar{u}(t,x)), \\ \hspace{4cm} t \in [0,T), x \in \mathbb{R}^d, \\ \bar{u}(T,x) = e^{\lambda T}\xi^{-1}(x)h(x), \qquad x \in \mathbb{R}^d, \end{cases} \tag{7}$$

where, for $i = \overline{1,k}$,
$$(\breve{\mathcal{E}}_t \bar{u})_i(t,x) = (\mathcal{L}_t \bar{u})_i(t,x) + \langle (\sigma\sigma^*)(t,x)\xi^{-1}(x)\nabla\xi(x), \nabla\bar{u}_i(t,x)\rangle$$

and

$$\begin{aligned}\bar{f}(t,x,\bar{u}) &= e^{\lambda t}\xi^{-1}(x)f(t,x,e^{-\lambda t}\xi(x)\bar{u}) - \lambda\bar{u} \\ &\quad + \tfrac{1}{2}\operatorname{Tr}[(\sigma\sigma^*)(t,x)\xi^{-1}(x)\mathrm{D}^2\xi(x)]\bar{u} + \langle b(t,x), \xi^{-1}(x)\nabla\xi(x)\rangle\bar{u}.\end{aligned}$$

Since $\lim_{|x|\to+\infty} \bar{u}(t,x) = \lim_{|x|\to+\infty} \bar{v}(t,x) = 0$ and $\bar{u}(T,x) = \bar{v}(T,x)$, there exists $(t_0, x_0) \in [0,T) \times \mathbb{R}^d$ such that

$$\theta := |\bar{u}(t_0,x_0) - \bar{v}(t_0,x_0)| \geq |\bar{u}(t,x) - \bar{v}(t,x)| \qquad \forall (t,x) \in [0,T] \times \mathbb{R}^d.$$

Let us suppose that $u \neq v$. This implies that $\theta > 0$. We set

$$z_0 = \frac{1}{\theta}(\bar{u}(t_0,x_0) - \bar{v}(t_0,x_0)).$$

Then,

$$\theta := \langle \bar{u}(t_0,x_0) - \bar{v}(t_0,x_0), z_0\rangle \geq \langle \bar{u}(t,x) - \bar{v}(t,x), z_0\rangle \qquad \forall (t,x) \in [0,T] \times \mathbb{R}^d.$$

We let, for $\alpha > 0$,

$$\Phi_\alpha(t,x,s,y) := \langle \bar{u}(t,x) - \bar{v}(s,y), z_0\rangle - \frac{\alpha}{2}(|t-s|^2 + |x-y|^2).$$

Since $\limsup_{|x|\vee|y|\to+\infty}\Phi_\alpha(t,x,s,y)\leq 0$ and $\Phi_\alpha(t_0,x_0,t_0,x_0)=\theta>0$, there exists $(t_\alpha,x_\alpha,s_\alpha,y_\alpha)\in([0,T]\times\mathbb{R}^d)^2$ such that

$$M_\alpha := \Phi_\alpha(t_\alpha,x_\alpha,s_\alpha,y_\alpha) = \sup\Phi_\alpha.$$

By Lemma 3.1 from [2], we have:

(a)
$$\lim_{\alpha\to+\infty}\alpha(|x_\alpha-y_\alpha|^2+|t_\alpha-s_\alpha|^2)=0;$$

(b) whenever $(\hat{t},\hat{x})$ is an accumulation point for $(t_\alpha,x_\alpha)$ as $\alpha\to+\infty$, we have that

$$\lim_{\alpha\to+\infty}M_\alpha = \langle\bar{u}(\hat{t},\hat{x})-\bar{v}(\hat{t},\hat{x}),z_0\rangle = \theta. \tag{8}$$

Since, for large $\alpha$, $(t_\alpha,x_\alpha)$ remains in a compact subset of $[0,T]\times\mathbb{R}^d$, there is at least one accumulation point $(\hat{t},\hat{x})$ for $(t_\alpha,x_\alpha)$ as $\alpha\to+\infty$. We can suppose, without loss of generality, that $(t_\alpha,x_\alpha)\to(\hat{t},\hat{x})$. Of course, from (8), $\hat{t}<T$, and we can also assume that $t_\alpha,s_\alpha\in[0,T)$ for every $\alpha$. Another consequence of (8) is that $\bar{u}(\hat{t},\hat{x})-\bar{v}(\hat{t},\hat{x})=\theta z_0$. Indeed,

$$|\bar{u}(\hat{t},\hat{x})-\bar{v}(\hat{t},\hat{x})|\leq\theta=\langle\bar{u}(\hat{t},\hat{x})-\bar{v}(\hat{t},\hat{x}),z_0\rangle\leq|\bar{u}(\hat{t},\hat{x})-\bar{v}(\hat{t},\hat{x})|,$$

from which we conclude that $\bar{u}(\hat{t},\hat{x})-\bar{v}(\hat{t},\hat{x})=|\bar{u}(\hat{t},\hat{x})-\bar{v}(\hat{t},\hat{x})|z_0=\theta z_0$. Let us now apply now Theorem 3.2 from [2], which asserts that, for every $\alpha>0$, there exist $X,Y\in S^d$ such that:

(a) $(\alpha(t_\alpha-s_\alpha),\alpha(x_\alpha-y_\alpha),X_\alpha)\in\bar{\mathcal{P}}^{2,+}\langle\bar{u},z_0\rangle(t_\alpha,x_\alpha);$

(b) $(\alpha(t_\alpha-s_\alpha),\alpha(x_\alpha-y_\alpha),Y_\alpha)\in\bar{\mathcal{P}}^{2,-}\langle\bar{v},z_0\rangle(s_\alpha,y_\alpha);$

(c) $\begin{pmatrix}X_a & 0\\ 0 & -Y_\alpha\end{pmatrix}\leq 3\alpha\begin{pmatrix}I & -I\\ -I & I\end{pmatrix}.$

From the definition of the viscosity solution and Remark 5, we have:

$$\alpha(t_\alpha-s_\alpha)+\bar{\mathcal{E}}_{\alpha(x_\alpha-y_\alpha),X_\alpha}(t_\alpha,x_\alpha)+\langle\bar{f}(t_\alpha,x_\alpha,\bar{u}(t_\alpha,x_\alpha)),z_0\rangle$$
$$\geq e^{\lambda t_\alpha}\xi^{-1}(x_\alpha)\varphi'_*(e^{-\lambda t_\alpha}\xi(x_\alpha)\bar{u}(t_\alpha,x_\alpha);z_0),$$
$$\alpha(t_\alpha-s_\alpha)+\bar{\mathcal{E}}_{\alpha(x_\alpha-y_\alpha),Y_\alpha}(s_\alpha,y_\alpha)+\langle\bar{f}(s_\alpha,y_\alpha,\bar{v}(s_\alpha,y_\alpha)),z_0\rangle$$
$$\leq e^{\lambda s_\alpha}\xi^{-1}(y_\alpha)\varphi',^*(e^{-\lambda s_\alpha}\xi(y_\alpha)\bar{v}(s_\alpha,y_\alpha);z_0),$$

where, for $q\in\mathbb{R}^d$, $X\in S^d$,

$$\bar{\mathcal{E}}_{q,X}(t,x)=\tfrac{1}{2}\mathrm{Tr}[(\sigma\sigma^*)(t,x)X]+\langle b(t,x),q\rangle+\langle(\sigma\sigma^*)(t,x)\xi^{-1}(x)\nabla\xi(x),q\rangle.$$

Subtracting the two inequalities, we obtain

$$\begin{aligned}
&\acute{\mathcal{E}}_{\alpha(x_\alpha-y_\alpha),X_\alpha}(t_\alpha,x_\alpha) - \acute{\mathcal{E}}_{\alpha(x_\alpha-y_\alpha),Y_\alpha}(s_\alpha,y_\alpha) \\
&\quad + \langle \bar{f}(t_\alpha,x_\alpha,\bar{u}(t_\alpha,x_\alpha)) - \bar{f}(s_\alpha,y_\alpha,\bar{v}(s_\alpha,y_\alpha)), z_0 \rangle \\
&\geq e^{\lambda t_\alpha} \xi^{-1}(x_\alpha) \varphi'_*(e^{-\lambda t_\alpha} \xi(x_\alpha) \bar{u}(t_\alpha,x_\alpha); z_0) \\
&\quad - e^{\lambda s_\alpha} \xi^{-1}(y_\alpha) \varphi'^{,*}(e^{-\lambda s_\alpha} \xi(y_\alpha) \bar{v}(s_\alpha,y_\alpha); z_0).
\end{aligned} \qquad (9)$$

From point (c) of Lemma 8, we get

$$\liminf_{\alpha \to +\infty} e^{\lambda t_\alpha} \xi^{-1}(x_\alpha) \varphi'_*(e^{-\lambda t_\alpha} \xi(x_\alpha) \bar{u}(t_\alpha,x_\alpha); z_0)$$

$$\geq e^{\lambda \hat{t}} \xi^{-1}(\hat{x}) \varphi'_*(e^{-\lambda \hat{t}} \xi(\hat{x}) \bar{u}(\hat{t},\hat{x}); z_0),$$

$$\limsup_{\alpha \to +\infty} e^{\lambda s_\alpha} \xi^{-1}(y_\alpha) \varphi'^{,*}(e^{-\lambda s_\alpha} \xi(y_\alpha) \bar{v}(s_\alpha,y_\alpha); z_0)$$

$$\leq e^{\lambda \hat{t}} \xi^{-1}(\hat{x}) \varphi'^{,*}(e^{-\lambda \hat{t}} \xi(\hat{x}) \bar{v}(\hat{t},\hat{x}); z_0).$$

Now, by point (d) of Lemma 8 and assumptions (A.5) and (A.6),

$$\varphi'_*(e^{-\lambda \hat{t}} \xi(\hat{x}) \bar{u}(\hat{t},\hat{x}); z_0) \geq \varphi'^{,*}(e^{-\lambda \hat{t}} \xi(\hat{x}) \bar{u}(\hat{t},\hat{x}); z_0)$$

since $\bar{u}(\hat{t},\hat{x}) - \bar{v}(\hat{t},\hat{x}) = \theta z_0$. It follows that

$$\begin{aligned}
\liminf_{\alpha \to +\infty} [&e^{\lambda t_\alpha} \xi^{-1}(x_\alpha) \varphi'_*(e^{-\lambda t_\alpha} \xi(x_\alpha) \bar{u}(t_\alpha,x_\alpha); z_0) \\
&\quad - e^{\lambda s_\alpha} \xi^{-1}(y_\alpha) \varphi'^{,*}(e^{-\lambda s_\alpha} \xi(y_\alpha) \bar{v}(s_\alpha,y_\alpha); z_0)] \\
&\geq e^{\lambda \hat{t}} \xi^{-1}(\hat{x}) \varphi'_*(e^{-\lambda \hat{t}} \xi(\hat{x}) \bar{u}(\hat{t},\hat{x}); z_0) \\
&\quad - e^{\lambda \hat{t}} \xi^{-1}(\hat{x}) \varphi'^{,*}(e^{-\lambda \hat{t}} \xi(\hat{x}) \bar{v}(\hat{t},\hat{x}); z_0) \geq 0.
\end{aligned}$$

On the other hand, we have that

$$\begin{aligned}
&\acute{\mathcal{E}}_{\alpha(x_\alpha-y_\alpha),X_\alpha}(t_\alpha,x_\alpha) - \acute{\mathcal{E}}_{\alpha(x_\alpha-y_\alpha),Y_\alpha}(s_\alpha,y_\alpha) \\
&= \tfrac{1}{2} \mathrm{Tr}[\sigma\sigma^*(t_\alpha,x_\alpha) X_\alpha - \sigma\sigma^*(s_\alpha,y_\alpha) Y_\alpha] \\
&\quad + \alpha \langle b(t_\alpha,x_\alpha) - b(s_\alpha,y_\alpha) \\
&\qquad + \sigma\sigma^*(t_\alpha,x_\alpha) \xi^{-1}(x_\alpha) \nabla \xi(x_\alpha) - \sigma\sigma^*(s_\alpha,y_\alpha) \xi^{-1}(y_\alpha) \nabla \xi(y_\alpha), x_\alpha - y_\alpha \rangle \\
&\leq \tfrac{3}{2} \alpha |\sigma(t_\alpha,x_\alpha) - \sigma(s_\alpha,y_\alpha)|^2 \\
&\quad + \alpha |x_\alpha - y_\alpha| (|b(t_\alpha,x_\alpha) - b(s_\alpha,y_\alpha)| \\
&\qquad + |\sigma\sigma^*(t_\alpha,x_\alpha) \xi^{-1}(x_\alpha) \nabla \xi(x_\alpha) - \sigma\sigma^*(s_\alpha,y_\alpha) \xi^{-1}(y_\alpha) \nabla \xi(y_\alpha)|) \\
&\leq 3\alpha (L_1^2 + L_1)(|x_\alpha - y_\alpha|^2 + |t_\alpha - s_\alpha|^2),
\end{aligned}$$

where $L_1$ is the Lipschitz constant of $(\sigma, b, \sigma\sigma^*\xi^{-1}\nabla\xi)$. Hence,

$$\limsup_{\alpha\to+\infty}[\bar{\mathcal{E}}'_{\alpha(x_\alpha-y_\alpha),X_\alpha}(t_\alpha,x_\alpha) - \bar{\mathcal{E}}'_{\alpha(x_\alpha-y_\alpha),Y_\alpha}(s_\alpha,y_\alpha)] \leq 0.$$

Finally, from assumption (A.2),

$$\lim_{\alpha\to+\infty}\langle\bar{f}(t_\alpha,x_\alpha,\bar{u}(t_\alpha,x_\alpha)) - \bar{f}(s_\alpha,y_\alpha,\bar{v}(s_\alpha,y_\alpha)), z_0\rangle$$
$$= \langle\bar{f}(\hat{t},\hat{x},\bar{u}(\hat{t},\hat{x})) - \bar{f}(\hat{t},\hat{x},\bar{v}(\hat{t},\hat{x})), z_0\rangle$$
$$= (-\lambda + \tfrac{1}{2}\mathrm{Tr}[\sigma\sigma^*(\hat{t},\hat{x})\xi^{-1}(\hat{x})\mathrm{D}^2\xi(\hat{x})] + \langle b(\hat{t},\hat{x}), \xi^{-1}(\hat{x})\nabla\xi(\hat{x})\rangle)\theta$$
$$\quad + \mathrm{e}^{\lambda\hat{t}}\xi^{-1}(\hat{x})\langle f(\hat{t},\hat{x},\mathrm{e}^{-\lambda\hat{t}}\xi(\hat{x})\bar{u}(\hat{t},\hat{x})) - f(\hat{t},\hat{x},\mathrm{e}^{-\lambda\hat{t}}\xi(\hat{x})\bar{v}(\hat{t},\hat{x})), z_0\rangle$$
$$\leq (-\lambda + \tfrac{1}{2}\mathrm{Tr}[\sigma\sigma^*(\hat{t},\hat{x})\xi^{-1}(\hat{x})\mathrm{D}^2\xi(\hat{x})] + \langle b(\hat{t},\hat{x}), \xi^{-1}(\hat{x})\nabla\xi(\hat{x})\rangle + \gamma)\theta.$$

By taking

$$\lambda > \gamma + \sup[\tfrac{1}{2}\mathrm{Tr}[\sigma\sigma^*\xi^{-1}\mathrm{D}^2\xi] + \langle b, \xi^{-1}\nabla\xi\rangle],$$

we get

$$\lim_{\alpha\to+\infty}\langle\bar{f}(t_\alpha,x_\alpha,\bar{u}(t_\alpha,x_\alpha)) - \bar{f}(s_\alpha,y_\alpha,\bar{v}(s_\alpha,y_\alpha)), z_0\rangle < 0.$$

Passing to the limit as $\alpha\to+\infty$ in (9), we obtain a contradiction. Hence, $\bar{u} = \bar{v}$ and so $u = v$.

□

## 4. Proof of the existence theorem

In this section, we prove Theorem 6. This is obtained via a stochastic approach. Using a certain backward stochastic variational inequality, we will obtain a generalized Feynman–Kac representation formula for the viscosity solution of (1).

Let $\{W_t\colon t\geq 0\}$ be a $d$-dimensional standard Brownian motion defined on some complete probability space $(\Omega, \mathcal{F}, \mathbb{P})$. We denote by $\{\mathcal{F}_t\colon t\geq 0\}$ the natural filtration generated by $\{W_t\colon t\geq 0\}$ and augmented by $\mathcal{N}$, the set of $\mathbb{P}$-null events of $\mathcal{F}$:

$$\mathcal{F}_t = \sigma\{W_r\colon 0\leq r\leq t\}\vee\mathcal{N}.$$

### 4.1. Backward stochastic variational inequality

We consider the following backward stochastic variational inequality

$$\begin{cases}\mathrm{d}Y_s + F(s,Y_s,Z_s)\,\mathrm{d}s \in \partial\varphi(Y_s)\,\mathrm{d}s + Z_s\,\mathrm{d}W_s, & s\in[0,T],\\ Y_T = \xi,\end{cases} \quad (10)$$

where:

(I) $F\colon \Omega \times [0,T] \times \mathbb{R}^k \times \mathbb{R}^{k\times d} \to \mathbb{R}^k$ satisfies, for some $\alpha \in \mathbb{R}$, $\beta \geq 0$,

$$
\begin{aligned}
&\text{(i)} && F(\cdot,\cdot,y,z) \text{ is a progressively measurable stochastic process (p.m.s.p.),} \\
&\text{(ii)} && y \longmapsto F(\omega,t,y,z) \text{ is continuous, a.s.,} \\
&\text{(iii)} && \langle y-y', F(t,y,z) - F(t,y',z) \rangle \leq \alpha |y-y'|^2, \text{ a.s.,} \\
&\text{(iv)} && |F(t,y,z) - F(t,y,z')| \leq L\|z-z'\|, \text{ a.s.,} \\
&\text{(v)} && \mathbb{E}\bigg(\int_0^T F_R^{\#}(t)\,\mathrm{d}t\bigg)^2 < \infty, \\
& && \forall R \geq 0, \text{ where } F_R^{\#}(t) := \sup\{|F(t,y,0)|\colon |y| \leq R\},
\end{aligned}
\qquad (11)
$$

for all $t \in [0,T]$, $y,y' \in \mathbb{R}^k$, $z,z' \in \mathbb{R}^{k\times d}$;
(II) $\varphi\colon \mathbb{R}^k \to (-\infty,+\infty]$ is proper, convex and l.s.c.;
(III) the terminal date $\xi \in L^2_{\mathcal{F}_T}(\Omega;\mathbb{R}^k)$ is such that

$$\mathbb{E}(\varphi(\xi)) < \infty. \qquad (12)$$

For the proof of the next theorem see [9].

**Theorem 10.** *Let the assumptions* (I)–(III) *be satisfied. There then exists a unique triple* $\{(Y_t,Z_t,U_t)\colon t \in [0,T]\}$ *of p.m.s.p. such that:*

(a) $\mathbb{E}\bigg[\sup\limits_{s\in[0,T]}|Y_s|^2 + \int_0^T (|Z_s|^2 + |U_s|^2)\,\mathrm{d}s\bigg] < +\infty;$

(b) $\mathbb{E}\int_0^T \varphi(Y_s)\,\mathrm{d}s < +\infty;$

(c) $(Y_t,U_t) \in \partial\varphi, \qquad \mathbb{P}\otimes \mathrm{d}t \text{ a.e.} \qquad \text{on } \Omega \times [0,T];$

(d) $Y_t + \int_t^T U_s\,\mathrm{d}s = \xi + \int_t^T F(s,Y_s,Z_s)\,\mathrm{d}s - \int_t^T Z_s\,\mathrm{d}W_s \qquad \text{for all } t \geq 0, \text{ a.s.}$

Let

$$\mathcal{F}_s^t = \sigma\{W_r - W_t\colon t \leq r \leq s\} \vee \mathcal{N}.$$

From assumption (A.1) of Section 2, it follows (see, for example, [6, 8]) that for each $(t,x) \in [0,T] \times \mathbb{R}^d$ there exists a unique continuous $\{\mathcal{F}_s^t\}$-p.m.s.p. $\{X_s^{t,x}\colon s \in [0,T]\}$ solving the stochastic differential equation

$$X_s^{t,x} = x + \int_t^{s\vee t} b(r,X_r^{t,x})\,\mathrm{d}r + \int_t^{s\vee t} \sigma(r,X_r^{t,x})\,\mathrm{d}W_r. \qquad (13)$$

***Remark.*** Note that, in this section, one can relax the assumption (A.1) by requiring only the continuity of $b$ and $\sigma$ in $t \in [0,T]$, instead of Lipschitz continuity.

The following proposition summarizes some well-known properties of the solutions of equation (13).

**Proposition 11.** *Under assumption* (A.1), *the solution* $\{X_s^{t,x}\colon s \in [0,T]\}$ *of equation (13) satisfies that for all $p \geq 2$, there exists some constant $C = C(p,T,L) > 0$ such that, for all $t,\tilde{t} \in [0,T]$, $x,\tilde{x} \in \mathbb{R}^d$:*

(i)    $\mathbb{E} \sup_{s \in [0,T]} |X_s^{t,x}|^p \leq C(1+|x|^p);$

(ii)    $\mathbb{E} \sup_{s \in [0,T]} |X_s^{t,x} - X_s^{\tilde{t},\tilde{x}}|^p \leq C(1+|x|^p+|\tilde{x}|^p)(|x-\tilde{x}|^p + |t-\tilde{t}|^{p/2}).$

From Theorem 10 with $\xi = h(X_T^{t,x})$, $F(\omega,s,y,z) = f(s, X_s^{t,x}(\omega), y)$, it follows, under assumptions (A.1)–(A.4), that for $(t,x) \in [0,T] \times \mathbb{R}^d$, there exists a unique triple of $\{\mathcal{F}_s^t\}$-p.m.s.p., $\{(Y_s^{t,x}, Z_s^{t,x}, U_s^{t,x})\colon s \in [0,T]\}$, solving the backward stochastic variational inequality

$$Y_s^{t,x} + \int_s^T U_r^{t,x}\, \mathrm{d}r = h(X_T^{t,x}) + \int_s^T f(r, X_r^{t,x}, Y_r^{t,x})\, \mathrm{d}r - \int_s^T Z_r^{t,x}\, \mathrm{d}W_r \quad (14)$$

for all $s \in [t,T]$ a.s.

with $(Y_s^{t,x}, U_s^{t,x}) \in \partial\varphi$, $\mathbb{P} \times \mathrm{d}s$ a.e. on $\Omega \times [t,T]$, and $Y_s^{t,x} = Y_t^{t,x}$, $Z_s^{t,x} = U_s^{t,x} = 0 \ \forall s \in [0,t]$.

Moreover, we have the following properties of the solution of (14).

**Proposition 12.** *Under the assumptions* (A.1)–(A.4), *we have that there exists some constant $C > 0$ such that for all $t,\tilde{t} \in [0,T]$, $x,\tilde{x} \in \mathbb{R}^d$:*

(i)    $\mathbb{E} \sup_{s \in [0,T]} |Y_s^{t,x}|^2 \leq C(1+|x|^2);$

(ii)    $\mathbb{E} \sup_{s \in [0,T]} |Y_s^{t,x} - Y_s^{\tilde{t},\tilde{x}}|^2$

$$\leq \mathbb{E}\bigg[|h(X_T^{t,x}) - h(X_T^{\tilde{t},\tilde{x}})|^2 + \int_0^T |\mathbf{1}_{[t,T]}(r) f(r, X_r^{t,x}, Y_r^{t,x}) - \mathbf{1}_{[\tilde{t},T]}(r) f(r, X_r^{\tilde{t},\tilde{x}}, Y_r^{t,x})|^2\, \mathrm{d}r\bigg].$$

We define

$$u(t,x) = Y_t^{t,x}, \qquad (t,x) \in [0,T] \times \mathbb{R}^d, \quad (15)$$

which is a determinist quantity since $Y_t^{t,x}$ is $\mathcal{F}_t^t \equiv \sigma\{\mathcal{N}\}$-measurable.

From the Markov property, we have that

$$u(s, X_s^{t,x}) = Y_s^{t,x}. \tag{16}$$

For the proof of the next proposition, see [11].

**Proposition 13.** *Under the assumptions* (A.1)–(A.4), *the function $u$ defined by (15) satisfies:*

$$\begin{aligned}&\text{(i)} && u(t,x) \in \text{Dom}(\varphi), \forall (t,x) \in [0,T] \times \mathbb{R}^d; \\ &\text{(ii)} && |u(t,x)| \leq C(1+|x|^p), \forall (t,x) \in [0,T] \times \mathbb{R}^d; \\ &\text{(iii)} && u \in C([0,T] \times \mathbb{R}^d; \mathbb{R}^k).\end{aligned} \tag{17}$$

## 4.2. Connection with variational inequalities for systems of partial differential equations

We prove that the function $u$ defined above is a viscosity solution, in the sense of Definition 4, for the multidimensional parabolic variational inequality (1).

Theorem 6 is a consequence of the following result.

**Theorem 14.** *Let the assumptions* (A.1)–(A.4) *be satisfied. The function $u$ defined by (15) is then a viscosity solution for the multidimensional parabolic variational inequality (1).*

**Proof.** Let $z \in \mathbb{R}^k$, $(t,x) \in [0,T] \times \mathbb{R}^d$ and $(p,q,X) \in \mathcal{P}^{2,+}\langle u,z\rangle(t,x)$. Let us define, for simplicity,

$$V(t,x,u,p,q,X) = p + \tfrac{1}{2}\text{Tr}((\sigma\sigma^*)(t,x)X) + \langle b(t,x), q\rangle + \langle f(t,x,u), z\rangle.$$

Suppose, contrary to our claim, that

$$V(t,x,u(t,x),p,q,X) < \varphi'_*(u(t,x);z).$$

It follows, by the continuity of $f$, $u$, $b$, $\sigma$, that there exist $\varepsilon > 0$, $\delta > 0$ such that for all $|s-t| \leq \delta$, $|y-x| \leq \delta$,

$$V(s,y,u(s,y),p+\varepsilon,q+(X+\varepsilon I)(y-x), X+\varepsilon I) < \varphi'_*(u(s,y);z). \tag{18}$$

Now, since $(p,q,X) \in \mathcal{P}^{2,+}\langle u,z\rangle(t,x)$, there exists $0 < \delta' \leq \delta$ such that

$$\langle u(s,y), z\rangle < \Psi(s,y) \tag{19}$$

for all $s \in [0,T]$, $s \neq t$, $y \in \mathbb{R}^d$, $y \neq x$ with $0 < s-t \leq \delta'$, $|y-x| \leq \delta'$, where

$$\Psi(s,y) = \langle u(t,x), z\rangle + (p+\varepsilon)(s-t) + \langle q, y-x\rangle + \tfrac{1}{2}\langle (X+\varepsilon I)(y-x), y-x\rangle.$$

Let

$$\tau = \inf\{s > t\colon |X_s^{t,x} - x| \geq \delta'\}.$$

If we let

$$(\bar{Y}_s^{t,x}, \bar{Z}_s^{t,x}) = (\langle Y_s^{t,x}, z\rangle, \langle Z_s^{t,x}, z\rangle), \qquad t \leq s,$$

then

$$\begin{cases} \bar{Y}_s^{t,x} = \langle u(\tau, X_\tau^{t,x}), z\rangle + \int_s^\tau (\langle f(r, X_r^{t,x}, u(r, X_r^{t,x})), z\rangle - \langle U_r^{t,x}, z\rangle)\,\mathrm{d}r - \int_s^\tau \bar{Z}_r^{t,x}\,\mathrm{d}W_r, \\ \langle U_s^{t,x}, z\rangle \in [\varphi'_-(Y_s^{t,x}; z), \varphi'_+(Y_s^{t,x}; z)]. \end{cases}$$

From Itô's formula, it follows that

$$(\hat{Y}_s^{t,x}, \hat{Z}_s^{t,x}) := (\Psi(s, X_s^{t,x}), (\nabla\Psi\sigma)(s, X_s^{t,x})), \qquad t \leq s \leq t + \delta',$$

satisfies

$$\hat{Y}_s^{t,x} = \Psi(\tau, X_\tau^{t,x}) - \int_s^\tau \left[\frac{\partial\Psi}{\partial r}(r, X_r^{t,x}) + \mathcal{L}_r\Psi(r, X_r^{t,x})\right]\mathrm{d}r - \int_s^\tau \hat{Z}_r^{t,x}\,\mathrm{d}W_r.$$

Let $(\tilde{Y}_s^{t,x}, \tilde{Z}_s^{t,x}) = (\hat{Y}_s^{t,x} - \bar{Y}_s^{t,x}, \hat{Z}_s^{t,x} - \bar{Z}_s^{t,x})$.

We have

$$\begin{aligned}\tilde{Y}_s^{t,x} &= [\Psi(\tau, X_\tau^{t,x}) - \langle u(\tau, X_\tau^{t,x}), z\rangle] \\ &\quad + \int_s^\tau \left[-\frac{\partial\Psi}{\partial r}(r, X_r^{t,x}) - \mathcal{L}_r\Psi(r, X_r^{t,x})\right. \\ &\quad\left. - \langle f(r, X_r^{t,x}, u(r, X_r^{t,x})), z\rangle + \langle U_r^{t,x}, z\rangle\right]\mathrm{d}r - \int_s^\tau \tilde{Z}_r^{t,x}\,\mathrm{d}W_r.\end{aligned} \qquad (20)$$

We note that, from Lemma 8, point (c),

$$\langle U_s^{t,x}, z\rangle \geq \varphi'_-(u(s, X_s^{t,x}); z) \geq \varphi'_*(u(s, X_s^{t,x}); z), \qquad \mathbb{P} \otimes \mathrm{d}t \text{ a.e.}$$

Moreover, the choice of $\delta'$ and $\tau$ implies that

$$\Psi(\tau, X_\tau^{t,x}) > \langle u(\tau, X_\tau^{t,x}), z\rangle.$$

From (18), it follows that

$$\begin{aligned}\varphi'_*(u(s, X_s^{t,x}); z) &> V(s, X_s^{t,x}, u(s, X_s^{t,x}), p+\varepsilon, q+(X+\varepsilon I)(X_s^{t,x}-x), X+\varepsilon I) \\ &= \frac{\partial\Psi}{\partial s}(s, X_s^{t,x}) + \mathcal{L}_s\Psi(s, X_s^{t,x}) + \langle f(s, X_s^{t,x}, u(s, X_s^{t,x})), z\rangle.\end{aligned}$$

These inequalities and equation (20) imply that $\tilde{Y}_t^{t,x} > 0$ or, equivalently,

$$\Psi(t, x) > \langle u(t, x), z\rangle,$$

which contradicts the definition of $\Psi$. Hence, we have

$$V(t, x, u(t,x), p, q, X) \geq \varphi'_*(u(t,x); z).$$

This proves that $u$ is a viscosity solution of (1). $\square$

## Acknowledgements

The work of L.M., A.R. and A.Z. was partially undertaken during visits to the Université de Provence and was supported by grants ID 395/2007 and CNCSIS 1156/2006.